\documentclass[11pt]{article}
\date{}
\makeatother
\usepackage{amssymb,amsmath}
\usepackage{amsmath}
\usepackage{color}
\usepackage[colorlinks,linkcolor=blue,citecolor=blue]{hyperref}
\usepackage{cite}
\usepackage{graphicx}
\usepackage{tikz}
\usepackage[applemac]{inputenc}
\usepackage{pgf}
\usepackage{mathtools}
\usepackage{psfrag}
\usepackage{xcolor}
\usepackage{pifont}
\usetikzlibrary{patterns}%ª≠Õº–Ë“™

\numberwithin{equation}{section}
 \newtheorem{thm}{Theorem}[section]
\newtheorem{lem}[thm]{Lemma}

\newtheorem{pro}[thm]{Proposition}

\newcommand{\RT}{{\mathbb{R}^3}}
\newcommand{\W}{W^{1,p}(\mathbb{R}^3)}
\newcommand{\D}{D^{1,q}(\mathbb{R}^3)}
\newcommand{\Q}{D^{1,q}}
\newcommand{\INT}{\int_{\mathbb{R}^3}}

\allowdisplaybreaks

\begin{document}

\title{The quasilinear Schr\"odinger--Poisson system}

\author{Yao Du$^{\rm a,b}$ \ \  \ \ \ Jiabao Su$^{\rm b\ast}$ \ \ \ \ \   Cong Wang$^{\rm c}\footnote{Corresponding author: Jiabao Su and Cong Wang.\  \hfill\break\indent \ \  E-mail addresses: 1052591976@qq.com(Y. Du).\ sujb@cnu.edu.cn(J. Su). wc252015@163.com(C. Wang).} $\\
{\small $^{\rm a}$School of Science, Xihua University,}\\
{\small Chengdu 610039, People's Republic of China}\\
{\small  $^{\rm b}$School of Mathematical Sciences,
Capital Normal University}\\ {\small Beijing 100048, People's Republic of China}\\
{ \small $^{\rm c}$Department of Mathematics,	Sichuan University}\\
 {  \small Chengdu 610064, People's Republic of China}
 }
\maketitle
 \begin{abstract}
This paper deals with the  $(p,q)$--Schr\"odinger--Poisson system, which is new and has never been considered in the literature. The uniqueness of solutions of the quasilinear Poisson equation is obtained via the Minty--Browder theorem.  The variational framework of the quasilinear system  is built and the nontrivial solutions of the system are obtained via   the mountain pass theorem.
 \end{abstract}
 \section{Introduction}
 In this paper  we consider the following Schr\"odinger--Poisson system
\begin{eqnarray}\label{101}
\left \{\begin{array}{ll}
\displaystyle -\Delta_p u+|u|^{p-2}u+\lambda\phi |u|^{s-2}u=|u|^{r-2}u,&\mathrm{in} \ \mathbb{R}^3,\\
\displaystyle -\Delta_q \phi = |u|^s,  &\mathrm{in}\ \mathbb{R}^3,\\
\end{array}
\right.
\end{eqnarray}
where $\Delta_i u=\hbox{div}(|\nabla u|^{i-2}\nabla u) (i=p,q)$, $\lambda>0$, $1<p<3$, $p<r<p^*$,  $q$ and $s$ satisfy
\begin{eqnarray}\label{2.6}
  \max \left\{1,\frac{3p}{5p-3}\right\}<q<3, \ \ \ \max\left\{1,\frac{(q^*-1)p}{q^*}\right\}<s<\frac{(q^*-1)p^*}{q^*},
\end{eqnarray}
$p^*=\frac{3p}{3-p}$ and $q^*=\frac{3q}{3-q}.$ We note that the first inequality in \eqref{2.6}
implies that $q$ satisfies $\frac{3p}{5p-3}<q<3$ for $1<p<\frac{3}{2}$ and $1<q<3$ for $\frac{3}{2}\leqslant p<3$.
 The system \eqref{101} is a quasilinear elliptic system coupled by a Schr\"odinger equation of $p$-Laplacian and a Poisson equation of $q$-Laplacian.
The $p$-Laplacian operator appears in nonlinear fluid dynamics, the range of $p$ is related to the velocity of flow and materials. For more information about the physical origin of the $p$-Laplacian, we refer to \cite{2018BGKT}.
 The Schr\"odinger--Poisson system comes from the quantum mechanics models\cite{1981BBL, 1992CL, 1981Lieb} and from the semiconductor theory \cite{1987Lions,1990MRS}.
  Motivating by the pioneering work of Benci and Fortunato\cite{1998BF},  the classical Schr\"odinger--Poisson system $(p=q=s=2)$ has been extensively studied via the variational methods. We refer to, for example, the works  \cite{2008AA, 2008AR, 2008Am, 2002BF, 2004DM-1,2004DM-2,  2006Ruiz} and references therein.

 In his famous work \cite{2006Ruiz}, Ruiz explored the effect of the nonlinear local term for the system \eqref{101} in the case that $p=q=s=2$, i.e. the coordinate point $(2,2,2)$ in the figure F1.
 In \cite{2021DSW-1} the system \eqref{101} has been considered for the case that $\frac{4}{3}<p<\frac{12}{5}$, $q=2$ and $s=2$, while in  \cite{2021DSW-2} the system \eqref{101} has been considered  for the case that  $1<p<3,$ $q=2$ and $s=p$, see the blue and the red line segment in the figure F1, respectively.

\begin{minipage}{\textwidth}\label{fig1}
\centering
\begin{tikzpicture}[scale=1]
\draw[->] (-1,0) -- (4,0);
\draw[->] (0,-1) -- (0,4);

\filldraw [gray] (2,2) circle (2pt);
\draw(2, 2) node[above] {$(2,2)$};

\draw[red,thick] (1,1)--(3,3);
\draw(1,1) node[below] {$(1,1)$};
\draw(3, 3) node[above] {$(3,3)$};

\draw[blue,thick] (4/3,2)--(12/5,2);
\draw (4/3,2) node[below left] {$(\frac{4}{3},2)$};
\draw (12/5,2) node[below right] {$(\frac{12}{5},2)$};

\draw (4, 0) node[below] {$p$};
\draw (0, 4) node[left] {$s$};

\end{tikzpicture}

{F1.\ \ The case $q=2$}
\end{minipage}\\

 As far as we know, for $q\neq2$, the system \eqref{101} has never   been studied in $\mathbb{R}^3$ before in the mathematical literature. A natural question arises: when $q\neq2$, whether the nontrivial solutions for the system \eqref{101} exists or not? What is the effect of the parameters $r$ and $\lambda$ on the existence of the nontrivial solutions for the system \eqref{101}? This is the main object of our paper.
 Motivated by the above mentioned works, especially our recent works \cite{2021DSW-1, 2021DSW-2}, this paper deals with the case $q\not=2$, or to be more precise, $1<p<3$, $q$ and $s$ satisfy \eqref{2.6}. The figure F2 shows the function image of $s$ and $p$ as $q=p$, the figure F3 gives a three-dimensional  visual sense without restriction $q=p$. From the figure F3, we see that the case  $q=2$ is included in the present paper. We remark that the  figure  F4(resp. F5) is the graph of the upper (resp. lower) bound of $s$, in order to create a better observation effect, the view angle is different from figure F3.

\begin{minipage}{\textwidth}\label{fig1}
\centering
\begin{tikzpicture}[scale=1]
\draw[->] (-1,0) -- (4,0);
\draw[->] (0,-1) -- (0,6.5);

\draw[blue,dashed] (6/5,1)--(3/2,1);
\draw[blue,dashed] (6/5,1)--(1.2,0);
\draw (1.2, 0) node[below] {$\frac{6}{5}$};
\draw[blue,dashed] (6/5,1)--(1.2,0);
\draw (1.2, 0) node[below] {$\frac{6}{5}$};

\draw[blue] (1.5,1.81)--(1.5,1);
\draw[blue,dashed] (1.5,1)--(1.5,0);
\draw (1.5, 0) node[below] {$\frac{3}{2}$};

\fill[pattern=north east lines, pattern color=red] (1.2,1)--(1.5,1)--(1.5,1.81)--(1.2,1)--cycle;

\draw[blue,dashed] (1.5,1)--(3,3);

\draw [thick,domain=1.2:2.5] plot(\x, {(\x)^2-0.44});

\draw[blue,dashed] (3,6.049647)--(3,0);
\draw (3, 0) node[below] {$3$};

\draw[blue,dashed] (1.2,1)--(0,1);
\draw (0,1) node[left] {$1$};

\draw[blue,dashed] (1.5,1.81)--(0,1.81);
\draw (0,1.81) node[left] {$2$};

\draw[blue,dashed] (3,3)--(0,3);
\draw (0,3) node[left] {$3$};
\draw (2,6/3) node[below right] {$\frac{4p-3}{3}$};
\draw (2,3.56) node[above left] {$\frac{4p-3}{3-p}$};

\fill[pattern=north west lines,dashed, pattern color=green] (1.5,1)--(3,3)--(3,5.81)--(2.58,5.81)--(1.5,1.81)--(1.5,1)--cycle;

\draw (4, 0) node[below] {$p$};
\draw (0, 6.5) node[left] {$s$};

\end{tikzpicture}

{F2.\ \ The case $p=q$}
\end{minipage}
\begin{minipage}{\textwidth}
\centering
\includegraphics[width=10cm,height=5.5cm]{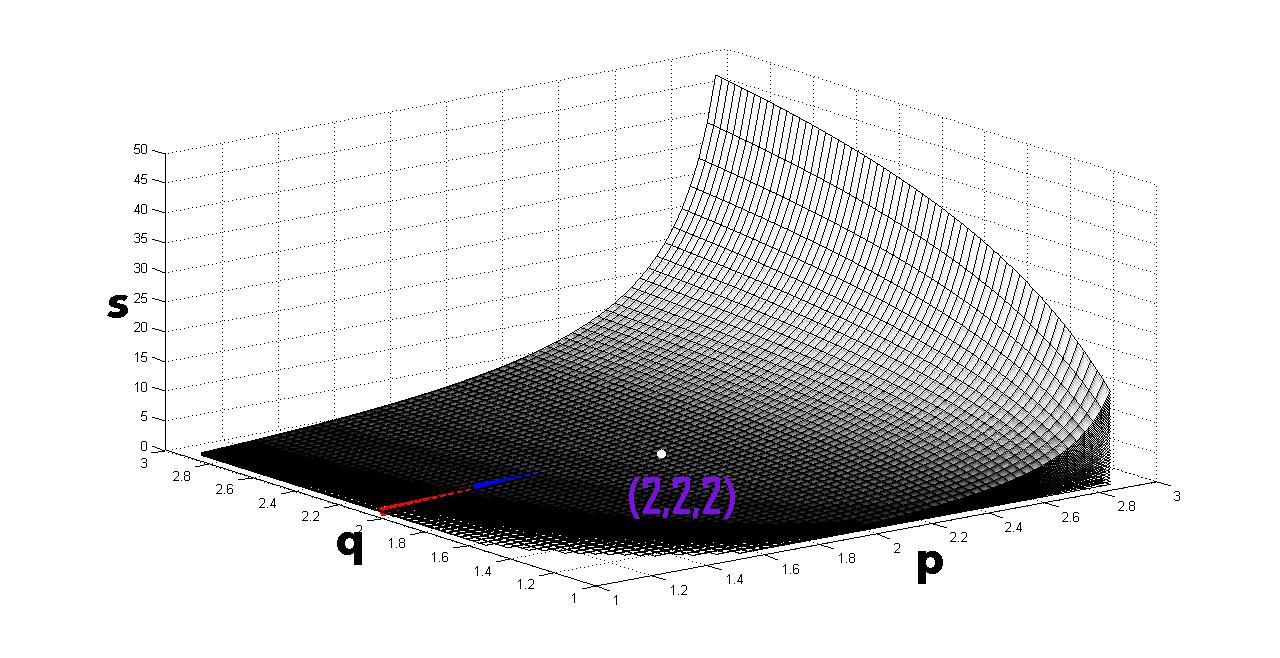}
\begin{center}
F3. The case involving $p\not =q$
\end{center}
\end{minipage}

\begin{figure}[htb]
\centering
\begin{minipage}[t]{0.3\textwidth}
\centering
\includegraphics[width=5.5cm,height=4cm]{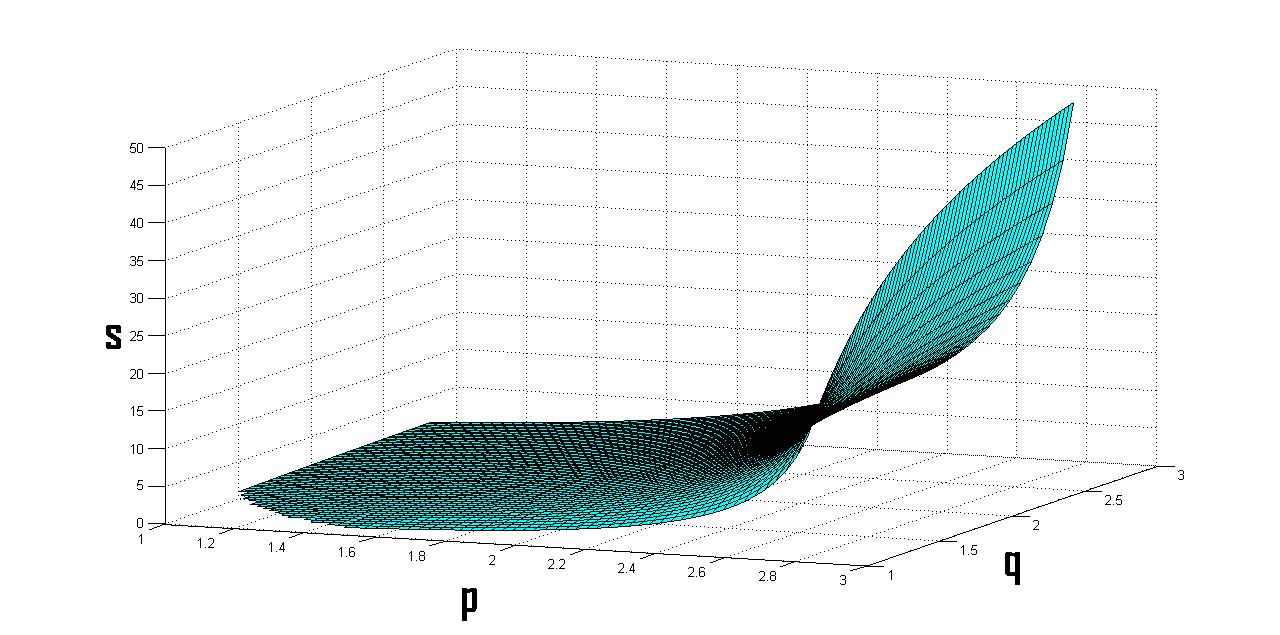}
\begin{center}
F4:\ $\frac{(q^*-1)p^*}{q^*}$.
\end{center}
\end{minipage}
\begin{minipage}[t]{0.3\textwidth}
\centering
\includegraphics[width=5.5cm,height=4cm]{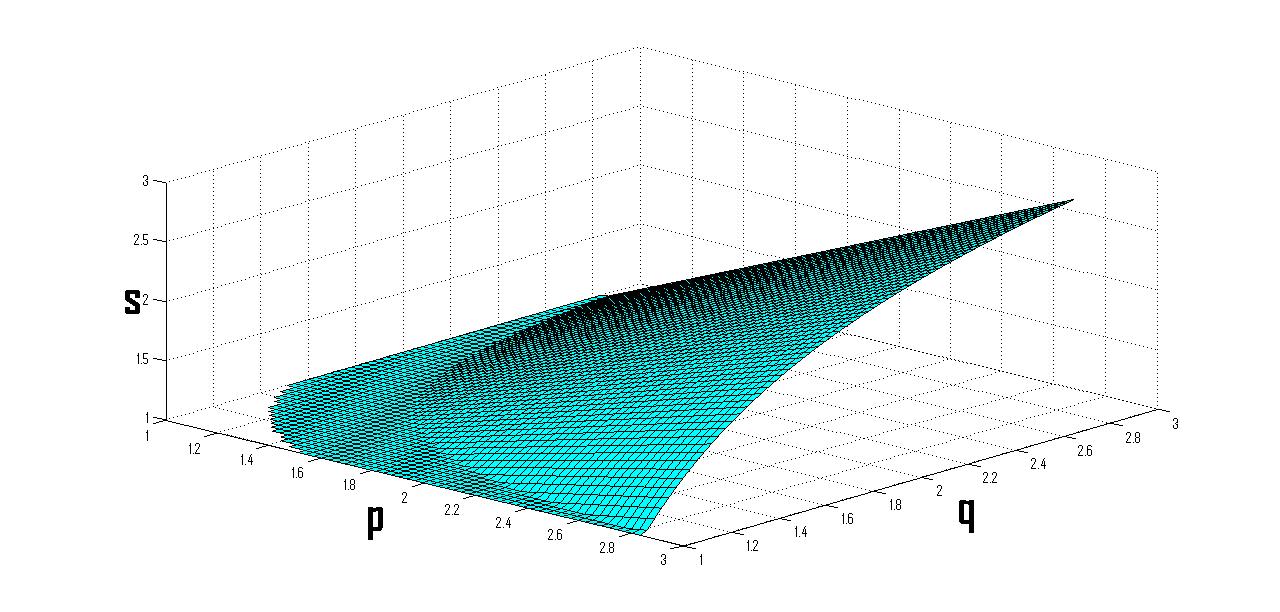}
\begin{center}
\ \ \ \ \ \ \ \ F5:\ $\max\left\{1,\frac{(q^*-1)p}{q^*}\right\}$.
\end{center}
\end{minipage}
\end{figure}

 The  stronger motivation than extending the numerical range stems from the nonlinearity of the
 operator $\Delta_q$ for  $q\not=2$ in the second equation of \eqref{101}.  This fact poses some essential obstacles in
  establishing the variational framework of \eqref{101}.

   In the case $q\ne 2$, there are two obstacles related to
   the second equation of \eqref{101}. The first one is that the Lax--Milgram theorem does not work in obtaining the unique solution $\phi_u\in \D$ for each $u\in \W$ in the second equation in \eqref{101}.  To overcome this difficulty, we apply the Minty--Browder theorem. The second one is that the unique solution
   $\phi_u$ may not have an explicit expression.  This prevents us from directly obtaining some further properties of $\phi_u$ including the non-negativeness
    and homogeneity, which are crucial in our discussion.  To overcome this obstacle, we rely on the uniqueness of the solution $\phi_u$ and the strict convexity of the energy functional for the second equation in \eqref{101}. Note that for $q=2$, $\phi_u$ can be written as
   $\phi_u(x)=\frac{1}{4\pi}\int_{\mathbb{R}^3}\frac{|u(y)|^s}{|x-y|}dy$, see \cite{2021DSW-2}.

   The real difficulty from $q\ne 2$ is the establishment of the variational framework of the system \eqref{101}. For $q=2$, the energy functional associated to \eqref{101}
   can be defined well via  by now almost a standard reduction procedure and it is of $C^1$, see e.g. \cite{1998BF, 2002BF, 2004DM-2, 2004DM-1, 2008AA, 2008AR, 2008Am, 2006Ruiz, 2021DSW-1, 2021DSW-2}.
   But for $q\ne 2$, it is a new task. We associate \eqref{101} with a functional $\mathcal{J}$ on $\W$  and  prove that $\mathcal{J}$ is of $C^1$ on $\W$.
    Then critical points of $\mathcal{J}$ correspond to solutions of \eqref{101}.
    The core is the $C^1$ regularity of the nonlinear operator $u\mapsto \phi_u$ plugged in $\mathcal{J}$.  The implicity of the function $\phi_u$ takes us a big difficulty. Inspired by \cite{2010YY}, we apply the basic methods in nonlinear functional analysis to  prove the  G\^{a}teaux differentiability of $\mathcal{J}$ and then to prove that the  G\^ateaux derivative of $\mathcal{J}$ is continuous.

As in the study of the classical Schr\"odinger--Poisson system, the difficulty appearing in the variational approach is the boundedness of the Palais--Smale sequences for $\mathcal{J}$. The key point is $r=\frac{q s}{q-1}$ and the delicate case is that $p<r<\frac{qs}{q-1}$.
  To achieve this property, we adopt the scaling technique in \cite{1997LJ} and the truncated technique in \cite{2006JS} to treat the case  $\frac{pq(1+s)}{p(q-1)+q}<r<\frac{qs}{q-1}$ and the case  $p<r\leqslant\frac{pq(1+s)}{p(q-1)+q}$, respectively.

 Hereafter we always assume that $1<p<3, q$ and $s$ satisfy \eqref{2.6}.  The main results are stated as follows.

 \begin{thm}\label{thm3.3}
 Let  $\frac{pq(1+s)}{p(q-1)+q}<r<p^*$.
 Then the system {\rm(\ref{101})} admits at least one nontrivial solution for any $\lambda>0$.
\end{thm}
 \begin{thm}\label{thm3.1}
 Let $p<r<p^*$.
 Then there exists $\lambda^*>0$ such that the system {\rm(\ref{101})} admits at least one nontrivial solution for any $0<\lambda<\lambda^*$.
\end{thm}

  To be honest, by now we do not know what would be the natural phenomena that
  the quasilinear system \eqref{101} could be used to explain.  We hope that this system would be useful to explore some known or unknown scientific phenomena or problems in the future.
  To the best of our knowledge, this is the first time that \eqref{101} is considered in $\mathbb{R}^3$ for $q\neq2$ in the mathematical literature.  We emphasize that the system \eqref{101} is a new model of variational methods and  it generalizes the system studied in {\rm\cite{2021DSW-1, 2021DSW-2, 2006Ruiz}}.

 We also note that the existence results of this paper coincide with those in {\rm\cite{2021DSW-1, 2021DSW-2, 2006Ruiz}} in the case $q=2$, $s=2$ or $s=p$. In particular, we extend the results in {\rm\cite{2021DSW-1, 2021DSW-2}} to $q\ne 2$.

This  paper is organized as follows. In Section \ref{sec2} we establish the variational setting of the system \eqref{101}.
Section \ref{sec3} is devoted to the proofs of Theorems \ref{thm3.3} and \ref{thm3.1}.

 \section{Variational framework} \label{sec2}
 In this section the variational structure of the system \eqref{101} will be established. We will work on the function spaces listed below.

$\bullet \ \W$ is the usual Sobolev space endowed with the norm
$$ \|u\|=\left(\int_{\mathbb{R}^3}|\nabla u|^p+|u|^p dx\right)^{\frac1p}.$$

$\bullet \ D^{1,q}(\mathbb{R}^3)$ denotes the completion of  $C_0^\infty(\mathbb{R}^3)$ with respect to the norm $$\|u\|_{D^{1,q}}=\left(\int_{\mathbb{R}^3} |\nabla u|^q dx\right)^{\frac1q}.$$

$\bullet \  L^s(\mathbb{R}^3)$, for $1 \leqslant s<\infty$, denotes the Lebesgue space with the norm $$\|u\|_{s}=\left(\int_{\mathbb{R}^3}|u|^s dx\right)^{\frac1s}.$$

 We use $X'$ to denote the dual space of a Banach space $X$ with the norm $\|\cdot\|_{X'}$ and use  $\langle\cdot,\cdot\rangle$ to denote the duality product between $X'$ and $X$.
 We use $ C$ to denote a positive constant that can change from line  to line.

 Based on the continuous  embeddings $\W \hookrightarrow  L^\ell (\mathbb{R}^3)$ for all $\ell\in [p, p^*]$
 and $\D \hookrightarrow L^{q^*}(\mathbb{R}^3)$,  the variational framework of \eqref{101} will be built by steps.
 Keep in mind the conditions  $1<q<3$ and $\frac{(q^*-1)p}{q^*}<s<\frac{(q^*-1)p^*}{q^*}$. We have
  \begin{pro}\label{Pro1} For any $u\in \W$, there exists a unique $\phi_u\in \D$ solving
 \begin{eqnarray}\label{6.8}
  \displaystyle -\Delta_q \phi  = |u|^s.
  \end{eqnarray}
\end{pro}
{\it Proof } \ For  $u\in \W$ and $v\in \D$,  H\"{o}lder and Sobolev inequalities imply that
\begin{eqnarray}\label{6.7}
 \left|\int_{\mathbb{R}^3} |u|^svdx\right|\leqslant
\left(\int_{\mathbb{R}^3}|u|^{\frac{q^*s}{q^*-1}} dx\right)^{\frac{q^*-1}{q^*}}\left(\int_{\mathbb{R}^3}|v|^{q^*}dx\right)^{\frac{1}{q^*}}
\leqslant C\| u\|^s \|v\|_{\Q}.
\end{eqnarray} It follows that for each fixed $u\in \W$, the linear functional
\begin{eqnarray*}
\mathcal{L}(v)=\int_{\mathbb{R}^3} |u|^svdx, \ \ \ v\in \D
 \end{eqnarray*} is well-defined and is continuous on $\D$.
 For $q=2$, by   the Lax--Milgram theorem, we immediately obtain the conclusion.
Notice that the Lax--Milgram theorem does not work for $q\neq2.$
 Instead, we apply the Minty--Browder theorem \cite{2011Brezis}.
 We   show that $\displaystyle -\Delta_q: \D \rightarrow (\D)'$ is a continuous map and satisfies
\begin{eqnarray}\label{2.4}
&& \langle-\Delta_q v_1-(-\Delta_q v_2),v_1-v_2\rangle>0, \ \ \ \forall \ v_1,v_2\in\D, \ v_1\neq v_2,
\end{eqnarray} \begin{eqnarray}
 \label{2.5}
\lim_{\|v\|_{\Q}\rightarrow\infty}\frac{\langle-\Delta_qv,v\rangle}{\|v\|_{\Q}}=\infty.
\end{eqnarray}
It is obvious that \eqref{2.5} holds since $q>1$.
 Let $v_n\rightarrow v$ in $\D$. Then $\{|\nabla v_n|^{q-2}\nabla v_n\}$ is bounded in $L^{\frac{q}{q-1}}(\RT)$
and $\nabla v_n(x)\rightarrow \nabla v(x)$  a.e. $x\in\RT$.  It follows from \cite[Proposition 5.4.7]{2013Willem} that
\begin{eqnarray*}
|\nabla v_n|^{q-2}\nabla v_n\rightharpoonup|\nabla v|^{q-2}\nabla v \ \ \mbox{in} \ \ L^{\frac{q}{q-1}}(\RT).
\end{eqnarray*}
Combining this with $\left\||\nabla v_n|^{q-2}\nabla v_n\right\|_{\frac{q}{q-1}}\rightarrow\left\||\nabla v|^{q-2}\nabla v\right\|_{\frac{q}{q-1}},$  we conclude that
\begin{eqnarray}\label{6.88}
\left(\INT \left||\nabla  v_n|^{q-2}\nabla v_n-|\nabla v|^{q-2}\nabla v\right|^{\frac{q}{q-1}}dx\right)^{\frac{q-1}{q}}\rightarrow0.
\end{eqnarray}
For any $\varphi\in\D$, by H\"{o}lder and Sobolev inequalities, we have
\begin{eqnarray}\label{7.1}
\left.
  \begin{array}{ll} \displaystyle
&{\langle-\Delta_q v_n-(-\Delta_q v),\varphi\rangle}\\[2mm]=&\displaystyle \INT (|\nabla v_n|^{q-2}\nabla v_n-|\nabla v|^{q-2}\nabla v)\nabla \varphi dx\\ \leqslant&\displaystyle\left(\INT\left||\nabla v_n|^{q-2}\nabla v_n-|\nabla v|^{q-2}\nabla v\right|^{\frac{q}{q-1}}dx\right)^{\frac{q-1}{q}}\|\varphi\|_{\Q}.
 \end{array}
\right.\label{7.1}
\end{eqnarray}
By \eqref{7.1} and \eqref{6.88}, we conclude that
\begin{eqnarray*}
&& \left\|-\Delta_q v_n-(-\Delta_q v)\right\|_{(\D)'}\\ &=&\sup\left\{|\langle-\Delta_q v_n-(-\Delta_q v),\varphi\rangle|: \ \varphi\in\D,\|\varphi\|_{\Q}=1\right\}\\
&\leqslant & \left(\INT \left||\nabla v_n|^{q-2}\nabla v_n-|\nabla v|^{q-2}\nabla v \right|^{\frac{q}{q-1}}dx\right)^{\frac{q-1}{q}}\rightarrow0.
\end{eqnarray*}
Therefore $-\Delta_q$ is continuous on $\D$.
We cite the elementary inequality: there exists $c_q>0$
  such that for any  $x, y \in \mathbb{R}^3$,  \begin{eqnarray} \label{252} \left\{\begin{array}{ll}
\langle|x|^{q-2}x-|y|^{q-2}y, x-y\rangle_{\mathbb{R}^3} \geqslant c_q |x-y|^q   &  \hbox{for} \ 2\leqslant q<3, \\
(|x|+|y|)^{2-q}\left\langle|x|^{q-2}x-|y|^{q-2}y, x-y\right\rangle_{\mathbb{R}^3}  \geqslant c_q |x-y|^2  & \hbox{for} \ 1<q<2,
  \end{array}
 \right.
 \end{eqnarray}
 where $\langle\cdot, \cdot\rangle_{\mathbb{R}^3}$ denotes the standard inner product in $\mathbb{R}^3$.
It follows from \eqref{252} that
\begin{eqnarray}\label{5.6}
c_q \|v_1 -v_2\|_{\Q}^{q}\leqslant \langle-\Delta_qv_1-(-\Delta_q v_2),v_1-v_2\rangle \ \ \mbox{for} \ 2\leqslant q<3,
 \end{eqnarray} and for $1<q<2$,
\begin{eqnarray}
\left.
  \begin{array}{ll} \displaystyle
 & c_q^{\frac{q}{2}} \|v_1 -v_2\|_{\Q}^{q} \\[3mm] \leqslant &  \displaystyle \int_{\RT} (T(v_1, v_2))^{\frac{q}{2}} \left(|\nabla v_1|+|\nabla v_2|\right)^{\frac{q(2-q)}{2}} dx\\
   \leqslant & \displaystyle
  \left(\langle-\Delta_q v_1-(-\Delta_q v_2), v_1-v_2\rangle\right)^{\frac{q}{2}} \left(\int_{\RT}\left(|\nabla v_1|+|\nabla v_2|\right)^{q} dx \right)^{\frac{2-q}{2}},
   \end{array}
\right.\label{5.7}
 \end{eqnarray}
 where
 \begin{eqnarray}\label{7.3}
T(v_1,v_2)=\left\langle|\nabla v_1|^{q-2}\nabla v_1-|\nabla v_2|^{q-2}\nabla v_2, \nabla v_1-\nabla v_2\right\rangle_{\mathbb{R}^3}.
 \end{eqnarray}
 From \eqref{5.6} and \eqref{5.7}, we deduce that \eqref{2.4} holds. By the continuity of the map $ -\Delta_q$, \eqref{2.4} and \eqref{2.5}, in view of $\mathcal{L}\in(\D)'$,   the conclusion follows from the Minty--Browder theorem.
The proof is complete.  \hfill$\Box$

 It is hard to give an explicit expression of $\phi_u$,  we can prove the following properties of $\phi_u$ via the uniqueness of the solution of \eqref{6.8}.
 \begin{pro}\label{Pro2} \ For $u\in \W,$ the solution $\phi_u$ of \eqref{6.8} has the following properties.
\begin{enumerate}
\item[{\rm (i) }]$\displaystyle\int_{\mathbb{R}^3} \left(\frac{1}{q} |\nabla \phi_u|^q - |u|^s\phi_u\right) dx
 =\min_{\phi\in\D} \int_{\mathbb{R}^3}\left(\frac{1}{q} |\nabla \phi|^q- |u|^s\phi \right) d x$, and
$\phi_u\geqslant0.$
 \item[{\rm (ii) }] For $t>0$, $\phi_{tu}=t^{\frac{s}{q-1}} \phi_u$ and $\phi_{u_t}(x)=t^{\frac{ks-q}{q-1}} \phi_u(tx)$, where $u_t(x)=t^ku(tx).$  Moreover $\phi_{u(\cdot+y)}=\phi_u(\cdot+y)$ for any
        $y\in\mathbb{R}^3$.
      \item[{\rm (iii) }]  $\|\phi_u\|_{\Q}\leqslant C\|u\|^{\frac{s}{q-1}}$, where C does not depend on $u.$
     \item[{\rm (iv) }] If $u_n\rightharpoonup u$ in $\W$, then $\phi_{u_n}\rightharpoonup \phi_u $ in $\D$ and $$\int_{\mathbb{R}^3}\phi_{u_n} |u_n|^{s-2}u_n\varphi dx\rightarrow \int_{\mathbb{R}^3}\phi_u|u|^{s-2}u\varphi dx, \ \ \forall \  \varphi\in\W.$$
    \item[{\rm (v) }] If $u_n\rightarrow u$ in $\W$, then $\phi_{u_n}\rightarrow\phi_u $ in $\D$.
 \end{enumerate}
\end{pro}
{\it Proof} \ (i) For each fixed $u\in\W$, we define the functional $I:\D\rightarrow\mathbb{R}$ as
 $$I(\phi)=\frac{1}{q}\int_{\mathbb{R}^3}|\nabla \phi|^qdx-\int_{\mathbb{R}^3}|u|^s\phi dx.$$
 It is easy to see that $I$ is of $C^1$ and $$I'(\phi)[v]=\int_{\mathbb{R}^3}|\nabla \phi|^{q-2}\nabla \phi\nabla vdx-\int_{\mathbb{R}^3}|u|^svdx.$$
By \eqref{2.4}, we have that
 \begin{eqnarray*}
(I'(w)-I'(v))[w-v]=\int_{\mathbb{R}^3}(|\nabla w|^{q-2}\nabla w-|\nabla v|^{q-2}\nabla v)(\nabla w-\nabla v)dx\geqslant0.
 \end{eqnarray*}
  Therefore $I$ is strictly convex.  Since $I$ is continuous and coercive, $I$ arrives   its global minimum
 uniquely at $\phi_u$ by Proposition \ref{Pro1}, i.e.,
$$\frac{1}{q}\int_{\mathbb{R}^3}|\nabla \phi_u|^qdx-\int_{\mathbb{R}^3}|u|^s\phi_udx=\min_{\phi\in\D}\left\{\frac{1}{q}\int_{\mathbb{R}^3}|\nabla \phi|^qdx-\int_{\mathbb{R}^3}|u|^s\phi dx\right\}.$$
Furthermore, it is easy to see that
$$\frac{1}{q}\int_{\mathbb{R}^3}|\nabla |\phi_u||^qdx-\int_{\mathbb{R}^3}|u|^s|\phi_u|dx\leqslant\frac{1}{q}\int_{\mathbb{R}^3}|\nabla \phi_u|^qdx-\int_{\mathbb{R}^3}|u|^s\phi_udx.$$
This shows that $|\phi_u|$ also achieves the minimum of $I$.
 Using again the uniqueness, $\phi_u=|\phi_u|\geqslant0.$

 (ii) For $t>0$, we have   \begin{eqnarray}\label{3.3}
 -\Delta_q(t\phi)=t^{q-1}(-\Delta_q \phi).
  \end{eqnarray}
  By \eqref{3.3}, we infer  that
 $$-\Delta_q\phi_{tu}=t^s|u|^s=t^s(-\Delta_q\phi_{u})=- \Delta_q(t^{\frac{s}{q-1}}\phi_{u}).$$
 Then it follows from the uniqueness that
 $$\phi_{tu}=t^{\frac{s}{q-1}}\phi_{u}.$$
 For $t>0$, we see that
    \begin{eqnarray}\label{3.6}
-\Delta_q\phi_{u_t}=|u_t|^s=t^{ks}|u(tx)|^s.
    \end{eqnarray}
 On the other hand, $-\Delta_q\phi_u=|u|^s$ implies that
   \begin{eqnarray}\label{3.4}
-\Delta_q\phi_u(tx)=|u(tx)|^s.
   \end{eqnarray}
An elementary computation gives that
   \begin{eqnarray}\label{3.5}
-\Delta_q\phi(tx)=t^{-q}(-\Delta_q(\phi(tx))).
   \end{eqnarray}
By \eqref{3.3}, \eqref{3.5} and \eqref{3.4}, we deduce that
  \begin{eqnarray}\label{3.8}
 -\Delta_q(t^{\frac{ks-q}{q-1}}\phi_u(tx))=t^{ks-q}(-\Delta_q(\phi_u(tx)))=t^{ks}(-\Delta_q\phi_u(tx))=t^{ks}|u(tx)|^s.
  \end{eqnarray}
 From \eqref{3.6} and \eqref{3.8}, using again the uniqueness, we infer that
 $$\phi_{u_t}(x)=t^{\frac{ks-q}{q-1}} \phi_u(tx).$$ In a similar way we have
 in view of the uniqueness that
 $$\phi_{u(\cdot+y)}(x)=\phi_u(x+y).$$

 (iii) \ Using H\"{o}lder and Sobolev inequalities, we conclude that
\begin{eqnarray*}
\|\phi_u\|_{\Q}^q=\int_{\mathbb{R}^3} |u|^s \phi_u dx\leqslant C\|u\|^s \| \phi_u\|_{\Q}.
\end{eqnarray*}

(iv) \ Let $u_n\rightharpoonup u$ in $\W$. Then we have that $u_n(x)\rightarrow u(x)$ a.e. $x\in\RT$
   and  $\{|u_n|^s\}$ is bounded in $L^{\frac{q^*}{q^*-1}}(\mathbb{R}^3)$.
Then, from \cite[Proposition 5.4.7]{2013Willem}, we deduce that
  \begin{eqnarray*}
   \int_{\mathbb{R}^3}(|u_n|^s-|u|^s)\varphi dx
 \rightarrow 0, \ \ \forall \  \varphi \in \D.
 \end{eqnarray*}
Combining this with \eqref{6.8}, we infer that $\phi_{u_n}\rightharpoonup \phi_u $ in $\D$. Then we have
\begin{eqnarray}\label{6.9}
\phi_{u_n}(x) u_n(x)\rightarrow \phi_u(x)u(x) \ \ \ \mbox{a.e.} \  x\in\RT.
   \end{eqnarray}
By H\"{o}lder and Sobolev inequalities, using (iii), we obtain that
  \begin{eqnarray*}
&& \int_{\mathbb{R}^3}\left|\phi_{u_n} |u_n|^{s-2}u_n\right|^{\frac{q^*s}{(s-1)q^*+1}}dx\\
&\leqslant&\left(\int_{\mathbb{R}^3}|\phi_{u_n}|^{q^*}dx\right)^{\frac{s}{(s-1)q^*+1}}
\left(\int_{\mathbb{R}^3}|u_n|^{\frac{q^*s}{q^*-1}}dx\right)^{\frac{(s-1)(q^*-1)}{(s-1)q^*+1}}\\
&\leqslant&C\|\phi_{u_n}\|_{\Q}^{\frac{q^*s}{(s-1)q^*+1}}\left(\int_{\mathbb{R}^3}|u_n|^{\frac{q^*s}{q^*-1}}dx\right)^{\frac{(s-1)(q^*-1)}{(s-1)q^*+1}}\\
&\leqslant&C\|{u_n}\|^{\frac{q^*s^2}{[(s-1)q^*+1](q-1)}}\left(\int_{\mathbb{R}^3}|u_n|^{\frac{q^*s}{q^*-1}}dx\right)^{\frac{(s-1)(q^*-1)}{(s-1)q^*+1}}.
\end{eqnarray*}
Thus $\{\phi_{u_n} |u_n|^{s-2}u_n\}$ is bounded in $L^{\frac{q^*s}{(s-1)q^*+1}}(\mathbb{R}^3)$. Together this with \eqref{6.9}, it follows from \cite[Proposition 5.4.7]{2013Willem} that
$$\int_{\mathbb{R}^3}\phi_{u_n} |u_n|^{s-2}u_n\varphi dx\rightarrow \int_{\mathbb{R}^3}\phi_u|u|^{s-2}u\varphi dx, \ \ \forall\ \varphi\in \W.$$

(v) \ Let $u_n\rightarrow u$ in $\W$. Then $|u_n|^s\rightarrow |u|^s$ in $L^{\frac{q^*}{q^*-1}}(\RT).$ From (iii) and (iv), we see that
 $\{\phi_{u_n}\}$ is bounded in $L^{q^*}(\mathbb{R}^3)$ and $\phi_{u_n}(x) \rightarrow \phi_u(x)\ \mbox{a.e.} \  {x\in\RT}$. Then from \cite[Proposition 5.4.7]{2013Willem}, we have
\begin{eqnarray} \label{6.1}
\int_{\mathbb{R}^3}(\phi_{u_n}-\phi_u)|u|^sdx\rightarrow0.
\end{eqnarray}
By H\"{o}lder inequality, using \eqref{6.1}, we infer that
\begin{eqnarray*}
\left|\int_{\mathbb{R}^3}|\nabla\phi_{u_n}|^q-|\nabla\phi_u|^qdx\right|
&=&\left|\int_{\mathbb{R}^3}\phi_{u_n}(|u_n|^s-|u|^s)+(\phi_{u_n}-\phi_u)|u|^sdx\right|\\
&\leqslant&\left\|\phi_{u_n}\right\|_{q^*}\left\||u_n|^s-|u|^s\right\|_{\frac{q^*}{q^*-1}}
+\left|\int_{\mathbb{R}^3}(\phi_{u_n}-\phi_u)|u|^sdx\right|\rightarrow0.
\end{eqnarray*}
Together this with (iv), we conclude that $\phi_{u_n}\rightarrow\phi_u $ in $\D$.
 The proof is complete.
 $\hfill\Box$

 By Proposition \ref{Pro2} and Sobolev inequality, the functional
 \begin{eqnarray*}
 \mathcal{J}(u)=\frac{1}{p}\int_{\mathbb{R}^3}\left(|\nabla u|^p+|u|^p\right) dx +\frac{\lambda(q-1)}{qs} \int_{\mathbb{R}^3}\phi_u|u|^sdx -\frac{1}{r}\int_{\mathbb{R}^3}|u|^{r}dx
 \end{eqnarray*} is well defined on $\W.$  We will prove that $\mathcal{J}\in C^1(\W, \mathbb{R})$ and
  critical points of $\mathcal{J}$ correspond to solutions of \eqref{101}.
  The key step is to prove that the operator $u \mapsto \phi_u \in C^1\left(\W,\D\right)$.
  For the case $q=2$, one may get the conclusion by the argument  in \cite{1998BF,2004DM-2} via the linearity of $(-\Delta)^{-1}$.  However, for the case $q\ne2$, the argument
  used in \cite{1998BF,2004DM-2} may not be effective to get
  the $C^1$ regularity of $u \mapsto \phi_u$.
  This causes that the argument for proving the $C^1$ regularity of   $\mathcal{J}$ is quite complicated.
  Inspired by \cite{2010YY}, we rely on the continuity of $u \mapsto \phi_u$ obtained in Proposition \ref{Pro2}(v).
\begin{pro}\label{lem2.4}
$\mathcal{J}\in C^1(\W,\mathbb{R})$ and for any $u, v\in\W$,
 \begin{eqnarray*}
\mathcal{J}'(u)[v]=\int_{\mathbb{R}^3}\left(|\nabla u|^{p-2}\nabla u\nabla v+|u|^{p-2}uv\right)dx +\lambda\int_{\mathbb{R}^3}\phi_u|u|^{s-2}uvdx-\int_{\mathbb{R}^3}|u|^{r-2}uvdx.
\end{eqnarray*}
\end{pro}
{\it Proof} \ We first prove that $\mathcal{J}$ is G\^{a}teaux differentiable. We will prove
 \begin{eqnarray}
\lim_{t\rightarrow0}\frac{\mathcal{J}(u+tv)-\mathcal{J}(u)}{t}=\mathcal{J}_G'(u)[v], \ \ \ u, v  \in\W,\label{4.9}
 \end{eqnarray}
where
$$\mathcal{J}_G'(u)[v]=\int_{\mathbb{R}^3}\left(|\nabla u|^{p-2}\nabla u\nabla v+|u|^{p-2}uv\right) dx +\lambda\int_{\mathbb{R}^3}\phi_u|u|^{s-2}uvdx-\int_{\mathbb{R}^3}|u|^{r-2}uvdx.$$
As a function of $v$, $\mathcal{J}_G'(u)[v]$ is a continuous linear functional on $\W$.
Set \begin{eqnarray*}
\mathcal{J}(u+tv)-\mathcal{J}(u)-t\mathcal{J}_G'(u)[v]=\mathcal{A}+\mathcal{B}+\lambda\mathcal{C},
\end{eqnarray*}
where
\begin{eqnarray*}
\mathcal{A}&=&\frac{1}{p} \|u+tv\|^p-\frac{1}{p}\|u\|^p-t\int_{\mathbb{R}^3}\left(|\nabla u|^{p-2}\nabla u\nabla v+|u|^{p-2} u v\right) dx,\\
\mathcal{B}&=&-\frac{1}{r}\int_{\mathbb{R}^3}|u+tv|^rdx+\frac{1}{r}\int_{\mathbb{R}^3}|u|^r dx+t\int_{\mathbb{R}^3}|u|^{r-2}uvdx,\\
\mathcal{C}&=&\frac{q-1}{qs} \int_{\mathbb{R}^3}\phi_{u+tv}|u+tv|^sdx-\frac{q-1}{qs} \int_{\mathbb{R}^3}\phi_u|u|^sdx-t\int_{\mathbb{R}^3}\phi_u |u|^{s-2}u vdx.
\end{eqnarray*}
It is easy to see that
\begin{eqnarray*}
\mathcal{A}=o(t), \ \ \ \mathcal{B}=o(t), \ \ \ t\to0.
\end{eqnarray*}
To prove \eqref{4.9}, it suffices to prove that
\begin{eqnarray}
\mathcal{C}=o(t), \ \ \ t\to0.\label{6.10}
\end{eqnarray}
From Proposition \ref{Pro2}(i), we see  that
\begin{eqnarray}\label{3.9}
\mathcal{M}(u,\phi_u)=\min_{\phi\in\D}\mathcal{M}(u,\phi),
\end{eqnarray}
where $$\mathcal{M}(u,\phi)=\frac{1}{q}\int_{\mathbb{R}^3}|\nabla \phi|^qdx-\int_{\mathbb{R}^3}|u|^s\phi dx.$$
Using $\mathcal{M}(u+tv,\phi_u)\geqslant \mathcal{M}(u+tv,\phi_{u+tv})$, we deduce  that
\begin{eqnarray}\label{4.11}
\mathcal{C}&=&-\frac{1}{s}\mathcal{M}(u+tv,\phi_{u+tv})-\frac{q-1}{qs} \int_{\mathbb{R}^3}\phi_u|u|^sdx-t\int_{\mathbb{R}^3}\phi_u|u|^{s-2}uvdx\nonumber\\
&\geqslant&-\frac{1}{s}\mathcal{M}(u+tv,\phi_u)-\frac{q-1}{q s} \int_{\mathbb{R}^3}\phi_u|u|^sdx-t\int_{\mathbb{R}^3}\phi_u|u|^{s-2}uvdx\\
&=&\frac{1}{s}\int_{\mathbb{R}^3}\phi_u \left(|u+tv|^s-|u|^s- ts |u|^{s-2}uv\right) dx.\nonumber
\end{eqnarray}
We claim that
\begin{eqnarray}
\int_{\mathbb{R}^3}\phi_u(|u+tv|^s-|u|^s-ts|u|^{s-2}uv)dx=o(t), \ \ t\to0.\label{4.10}
\end{eqnarray}
Indeed, in view of $s>1$, we have that for a.e. $x\in\RT$,
\begin{eqnarray}
\lim_{t\rightarrow0}\frac{\phi_u(x)(|u(x)+tv(x)|^s-|u(x)|^s-ts|u(x)|^{s-2}u(x)v(x))}{t}=0.\label{3.18}
\end{eqnarray}
We cite the elementary inequality: for any $q>0,$ there exists $C_q>0$ such that
\begin{eqnarray}  \label{3.10}
|a+b|^q\leqslant C_q(|a|^q+|b|^q),\ \ \forall \ a,b\in\mathbb{R}.
\end{eqnarray}
Then, as $s>1$, by the Lagrange theorem there exists $\theta\in \mathbb{R}$ such that $|\theta|\leqslant|t|$ and
\begin{eqnarray}\label{3.17}
&& \frac1t \left| \phi_u(x)[|u(x)+tv(x)|^s-|u(x)|^s-t s|u(x)|^{s-2}u(x)v(x)]\right|\nonumber\\
&\leqslant& s\phi_u(x)\left(\left||u(x)+\theta v(x)|^{s-2}(u(x)+\theta v(x))v(x)\right|+\left||u(x)|^{s-2}u(x)v(x)\right|\right)\\
&\leqslant& C\phi_u(x)\left(|u(x)|^{s-1}|v(x)|+|v(x)|^s\right).\nonumber
\end{eqnarray}
By H\"{o}lder and Sobolev  inequalities, using \eqref{3.10} and Proposition \ref{Pro2}(iii), we infer that
\begin{eqnarray}
\left.
  \begin{array}{ll}& \displaystyle
\left|\int_{\mathbb{R}^3}\phi_u(|u|^{s-1}|v|+|v|^s)dx\right|\\  \leqslant & \displaystyle \|\phi_u\|_{q^*}
\left(\int_{\mathbb{R}^3}(|u|^{s-1}|v|+|v|^s)^{\frac{q^*}{q^*-1}} dx\right)^{\frac{q^*-1}{q^*}}\\
 \leqslant & \displaystyle C\|\phi_u\|_{\Q}
\left(\int_{\mathbb{R}^3}(|u|^{s-1}|v|)^{\frac{q^*}{q^*-1}}+|v|^{\frac{q^*s}{q^*-1}} dx\right)^{\frac{q^*-1}{q^*}}\\
 \leqslant & \displaystyle C\|u\|^{\frac{s}{q-1}}
\left(\int_{\mathbb{R}^3}(|u|^{s-1}|v|)^{\frac{q^*}{q^*-1}}+|v|^{\frac{q^*s}{q^*-1}} dx\right)^{\frac{q^*-1}{q^*}}.
 \end{array}
\right.\label{3.16}
\end{eqnarray}
As $\max\left\{1,\frac{(q^*-1)p}{q^*}\right\}<s<\frac{(q^*-1)p^*}{q^*}$, by H\"{o}lder inequality and Sobolev embedding,
we have
\begin{eqnarray}\label{3.13}
\int_{\mathbb{R}^3}\left||u|^{s-1}v\right|^{\frac{q^*}{q^*-1}} dx \leqslant C \|u\|^{\frac{q^*(s-1)}{q^*-1}}\|v\|^{\frac{q^*}{q^*-1}}.
\end{eqnarray}
From \eqref{3.16} and \eqref{3.13}, the function $\phi_u(|u|^{s-1}|v|+|v|^s)$ is in $L^1(\RT).$ Combining this with \eqref{3.17} and \eqref{3.18}, by the dominated
convergence theorem, \eqref{4.10} follows. From \eqref{4.11} and \eqref{4.10}, we conclude
that
\begin{eqnarray}\label{4.6}
\mathcal{C}\geqslant o(t).
\end{eqnarray}
We now verify that
\begin{eqnarray}\label{4.5}
\mathcal{C}\leqslant o(t).
\end{eqnarray}
Actually, using $\mathcal{M}(u,\phi_u)\leqslant \mathcal{M}(u,\phi_{u+tv})$, by H\"{o}lder inequality, we deduce that
\begin{eqnarray}\label{3.15}
\mathcal{C}&=&\frac{q-1}{q s} \int_{\mathbb{R}^3}\phi_{u+tv}|u+tv|^sdx+\frac{1}{s}\mathcal{M}(u,\phi_u)-t\int_{\mathbb{R}^3}\phi_u|u|^{s-2}uvdx\nonumber\\
&\leqslant&\frac{q-1}{q s} \int_{\mathbb{R}^3}\phi_{u+tv}|u+tv|^sdx+\frac{1}{s}\mathcal{M}(u,\phi_{u+tv})-t\int_{\mathbb{R}^3}\phi_u|u|^{s-2}uvdx\nonumber\\
&=&\frac{1}{s}\int_{\mathbb{R}^3}[\phi_{u+tv}(|u+tv|^s-|u|^s-st|u|^{s-2}uv)
+st(\phi_{u+tv}-\phi_u)|u|^{s-2}uv]dx\\
&\leqslant&\frac{1}{s}\left(\int_{\mathbb{R}^3}|\phi_{u+tv}|^{q^*}dx\right)^{\frac{1}{q^*}}
\left(\int_{\mathbb{R}^3}\left||u+tv|^s-|u|^s-st|u|^{s-2}uv\right|^{\frac{q^*}{q^*-1}} dx\right)^{\frac{q^*-1}{q^*}}\nonumber\\
&&+t\left(\int_{\mathbb{R}^3}|\phi_{u+tv}-\phi_u|^{q^*}dx\right)^{\frac{1}{q^*}}
\left(\int_{\mathbb{R}^3}\left||u|^{s-2}uv\right|^{\frac{q^*}{q^*-1}} dx\right)^{\frac{q^*-1}{q^*}}.\nonumber
\end{eqnarray}
We first check that
\begin{eqnarray} \label{3.11}
\lim_{t\rightarrow0} \frac1t \displaystyle
\left(\int_{\mathbb{R}^3}\left||u+tv|^s-|u|^s-st|u|^{s-2}uv\right|^{\frac{q^*}{q^*-1}} dx\right)^{\frac{q^*-1}{q^*}} =0.
\end{eqnarray}
To prove \eqref{3.11}, it suffices to prove that
\begin{eqnarray} \label{3.14}
\lim_{t\rightarrow0}  t^{-\frac{q^*}{q^*-1}} \displaystyle
\int_{\mathbb{R}^3}
 \left||u+tv|^s-|u|^s-st|u|^{s-2}uv\right|^{\frac{q^*}{q^*-1}} dx =0.
\end{eqnarray}
Indeed, since $s>1$, it is clear that for a.e. $x\in\RT$,
\begin{eqnarray}\label{3.12}
\lim_{t\rightarrow0} \frac{|u(x)+tv(x)|^s-|u(x)|^s-st|u(x)|^{s-2}u(x)v(x)}{t}
 =0.
\end{eqnarray}
As in \eqref{3.17}, by the Lagrange theorem there exists $\theta\in \mathbb{R}$ such that $|\theta|\leqslant|t|$ and
\begin{eqnarray}\label{4.3}
\frac1t \left| |u(x)+tv(x)|^s-|u(x)|^s-st|u(x)|^{s-2}u(x)v(x) \right|
\leqslant C(|u(x)|^{s-1}|v(x)|+|v(x)|^s).
\end{eqnarray}
From \eqref{3.13}, we infer that
the function $(|u|^{s-1}|v|+|v|^s)^{\frac{q^*}{q^*-1}}$ is in $L^1(\RT).$
 Combining this with \eqref{3.12} and \eqref{4.3}, by the dominated
convergence theorem, \eqref{3.14}   follows.
Applying Proposition \ref{Pro2}(v) and Sobolev embedding theorem  we can deduce that
\begin{eqnarray}
 \lim_{t\rightarrow0}
 \int_{\mathbb{R}^3}|\phi_{u+tv}-\phi_u|^{q^*}dx =0.\label{4.2}
\end{eqnarray}
By \eqref{3.11} and \eqref{4.2}, we deduce from \eqref{3.15} that \eqref{4.5} holds. It follows from  \eqref{4.5} and \eqref{4.6} that \eqref{6.10} holds.

To end the proof, we need to prove the continuity of $\mathcal{J}_G': \W \rightarrow (\W)'.$  We only need to prove
that
 $$\Psi_G'(u)[v]  := \int_{\mathbb{R}^3}   \phi_u|u|^{s-2}u v dx $$  is continuous on $\W$.  Let $u_n\rightarrow u$ in $\W$. By  Proposition \ref{Pro2}(v) and Sobolev embedding theorem, we infer that $\phi_{u_n}\rightarrow\phi_u $ in $L^{q^*}(\RT)$. Moreover, $|u_n|^{s-2}u_n\rightarrow|u|^{s-2}u$ in $L^{\frac{q^*s}{(s-1)(q^*-1)}}(\RT)$.
Then, by H\"{o}lder inequality, we obtain that for any $v\in\W$,
\begin{eqnarray*}
&&{\langle\Psi_G'(u_n)-\Psi_G'(u),v\rangle}\nonumber\\
&=&\int_{\mathbb{R}^3}(\phi_{u_n}|u_n|^{s-2}u_n-\phi_u|u|^{s-2}u)vdx\nonumber\\
&\leqslant&
\int_{\mathbb{R}^3}|\phi_{u_n}(|u_n|^{s-2}u_n-|u|^{s-2}u)v|dx+
\int_{\mathbb{R}^3}|(\phi_{u_n}-\phi_u)|u|^{s-2}uv|dx\\
&\leqslant&
\left\||u_n|^{s-2}u_n-|u|^{s-2}u\right\|_{\frac{q^*s}{(s-1)(q^*-1)}}
\|\phi_{u_n}\|_{q^*}\|v\|_{\frac{q^*s}{q^*-1}}+\left\|\phi_{u_n}-\phi_u\right\|_{q^*}
\|u\|_{\frac{q^*s}{q^*-1}}^{s-1}
\|v\|_{\frac{q^*s}{q^*-1}}
\nonumber\\
&\rightarrow&0.\nonumber
\end{eqnarray*}
Therefore
\begin{eqnarray*}
&& \|\Psi_G'(u_n)-\Psi_G'(u)\|_{(\W)'}\\ &=&\sup\left\{|\langle\Psi_G'(u_n)-\Psi_G'(u),v\rangle|: \ v\in\W,\|v\|=1\right\}\\
&\leqslant& C \left\||u_n|^{s-2}u_n-|u|^{s-2}u\right\|_{\frac{q^*s}{(s-1)(q^*-1)}}\|\phi_{u_n}\|_{q^*}
 +C\|\phi_{u_n}-\phi_u\|_{q^*}
\|u\|_{\frac{q^*s}{q^*-1}}^{s-1}\\
&\rightarrow&0.
\end{eqnarray*}
 The proof is complete.\hfill $\Box$

  \begin{pro}\label{lem2.1}
Let $(u,\phi)\in\W\times\D$. Then $(u,\phi)$ is a solution of \eqref{101}
 if and only if  $u$ is a critical point of $\mathcal{J}$ and $\phi=\phi_u$.
     \end{pro}
{\it Proof} \ We first define the functional $\mathcal{H}: \W\times\D\rightarrow\mathbb{R}$ as
  $$\mathcal{H}(u,\phi)=\frac{1}{p}\int_{\mathbb{R}^3}\left(|\nabla u|^p+|u|^p\right) dx +\frac{\lambda}{s} \int_{\mathbb{R}^3}\phi|u|^sdx -\frac{\lambda}{qs} \int_{\mathbb{R}^3}|\nabla \phi|^qdx -\frac{1}{r}\int_{\mathbb{R}^3}|u|^{r}dx .$$
Observe that for any $(v,w)\in\W\times\D$,
 $$\partial_u \mathcal{H}(u,\phi)[v]=\int_{\mathbb{R}^3}\left(|\nabla u|^{p-2}\nabla u\nabla v+|u|^{p-2}uv\right) dx +\lambda\int_{\mathbb{R}^3}\phi|u|^{s-2}uvdx-\int_{\mathbb{R}^3}|u|^{q-2}uvdx,$$
  $$\partial_\phi \mathcal{H}(u,\phi)[w]=\frac{\lambda}{s}\int_{\mathbb{R}^3}|u|^sw dx-\frac{\lambda}{s} \int_{\mathbb{R}^3}|\nabla \phi|^{q-2}\nabla \phi \nabla w dx.$$
  Then we deduce that
   \begin{eqnarray*}
  (u, \phi)\  \mbox{is\ a\ solution\ of}\ \eqref{101}\ \Leftrightarrow \
\partial_u \mathcal{H} (u,\phi)=0 \
\mbox{and}\ \partial_\phi \mathcal{H}(u,\phi)=0.
   \end{eqnarray*}
 From the definition of $\mathcal{J}'(u)$ given in Proposition  \ref{lem2.4}, we infer that
   \begin{eqnarray*}
\partial_u \mathcal{H}(u,\phi)=0 \ \mbox{and}\ \partial_\phi \mathcal{H}(u,\phi)=0\ \
\Leftrightarrow \ \ \mathcal{J}'(u)=0 \ \mbox{and}\ \phi=\phi_u.
   \end{eqnarray*}
The proof is complete. \hfill $\Box$

 \section{Proofs of Theorems \ref{thm3.3} and \ref{thm3.1}} \label{sec3}
In this section we prove Theorems \ref{thm3.3} and \ref{thm3.1}. We first give a few of important lemmas.
\begin{lem}\label{lem1}
 Let $\{u_n\}\subset\W$ be a bounded  sequence such that $\mathcal{J}'(u_n)\rightarrow0$.
 Then there exists  $u\in \W$ so that, up to a subsequence, $\nabla u_n(x)\to\nabla u(x)$ a.e. $x\in\RT$.
\end{lem}
{\it Proof} \ Since $\{u_n\}$ is bounded in $\W$, there  exists  $u\in \W$ so that, up to a subsequence,
\begin{eqnarray}
  \left\{\begin{array}{ll}
u_n\rightharpoonup u &  \mbox{in} \ \W;\\
u_n\rightarrow u &  \mbox{in} \ L_{\rm loc}^l({\mathbb{R}^3})  \ \ \forall \ p\leqslant l<p^*;\\
u_n(x)\rightarrow u(x) &  \mbox{a.e.\ } \  {x\in\mathbb{R}^3}.
  \end{array}
\right.\label{7.2}
\end{eqnarray}
We follow some ideas in \cite{1992BM}. Let $\eta \in C_{0}^\infty(\mathbb{R}^3, [0, 1])$ satisfy $\eta|_{B_{R}}=1$ and $\hbox{supp}\ \eta \subset B_{2R}$, where  $B_R=\left\{x\in \mathbb{R}^3: |x|\leqslant R\right\}$.
From $\mathcal{J}'(u)\in(\W)'$, $(u_n-u)\eta\rightharpoonup0$ in $\W$ and $\mathcal{J}'(u_n)\rightarrow0$, we infer that
  \begin{eqnarray} \label{307}
  (\mathcal{J}'(u_n)-\mathcal{J}'(u))[(u_n-u)\eta]\to 0.
 \end{eqnarray}
By H\"older inequality, using \eqref{7.2}, we conclude that
\begin{eqnarray}\label{4.14}
&&\int_{\mathbb{R}^3} \left(|u_n|^{l-2}u_n-|u|^{l-2} u\right)(u_n-u) \eta dx=o(1),\ \ \forall \ p\leqslant l<p^*,\\
&&\int_{\mathbb{R}^3} |\nabla u_n|^{p-2}\nabla u_n \nabla \eta \ (u_n-u) dx=o(1).
\end{eqnarray}
Moreover, we deduce that
\begin{eqnarray}\label{4.17}
\INT(\phi_{u_n} |u_n|^{s-2}u_n-\phi_{u} |u|^{s-2}u) (u_n-u)\eta dx=o(1).
\end{eqnarray}
Indeed, by H\"{o}lder and Sobolev inequalities, using Proposition \ref{Pro2}(iii), we have
\begin{eqnarray*}
&&\INT(\phi_{u_n} |u_n|^{s-2}u_n-\phi_{u} |u|^{s-2}u) (u_n-u)\eta dx\nonumber\\
&\leqslant&\|\phi_{u_n}\|_{q^*}\left(\int_{B_{2R}}\left||u_n|^{s-2}u_n (u_n-u)\eta\right|^{\frac{q^*}{q^*-1}} dx\right)^{\frac{q^*-1}{q^*}}\nonumber\\
&&+\|\phi_u\|_{q^*}\left(\int_{B_{2R}}\left||u|^{s-2}u (u_n-u)\eta\right|^{\frac{q^*}{q^*-1}} dx\right)^{\frac{q^*-1}{q^*}}\nonumber\\
&\leqslant&C\|\phi_{u_n}\|_{\Q}\left(\int_{B_{2R}}\left||u_n|^{s-2}u_n (u_n-u)\eta\right|^{\frac{q^*}{q^*-1}} dx\right)^{\frac{q^*-1}{q^*}}\\
&&+C\|\phi_u\|_{\Q}\left(\int_{B_{2R}}\left||u|^{s-2}u (u_n-u)\eta\right|^{\frac{q^*}{q^*-1}} dx\right)^{\frac{q^*-1}{q^*}}\nonumber\\
&\leqslant&C(\|u_n\|^{\frac{s}{q-1}}\|u_n\|_{\frac{q^*s}{q^*-1}}^{s-1}+\|u\|^{\frac{s}{q-1}}\|u\|_{\frac{q^*s}{q^*-1}}^{s-1})
\left(\int_{B_{2R}}|u_n-u|^{\frac{q^*s}{q^*-1}} dx\right)^{\frac{q^*-1}{q^*s}}.
\end{eqnarray*}
Combining the above inequality and \eqref{7.2} we get \eqref{4.17}.
By  \eqref{307}--\eqref{4.17}, we infer that
\begin{eqnarray*}
\int_{\mathbb{R}^3}\left(|\nabla u_n|^{p-2}\nabla u_n-|\nabla u|^{p-2}\nabla u\right)\nabla (u_n-u)\eta dx=o(1).
\end{eqnarray*}
By \eqref{252}, $T(u_n,u)\geqslant0$ easily follows, where $T(\cdot,\cdot)$ is given in \eqref{7.3}.
Then, we see that
  \begin{eqnarray}\label{207}
\lim_{n\to\infty} \int_{B_{R}}T(u_n,u)dx =0.
\end{eqnarray}
Using again \eqref{252}, proceeding as in \eqref{5.6} and \eqref{5.7}, we deduce from \eqref{207}
that
 \begin{eqnarray*}
 \lim_{n\to\infty} \int_{B_R}  |\nabla u_n-\nabla u|^{p} dx  =0. \end{eqnarray*}
 Then, up to a subsequence, we get
\begin{eqnarray*}
\nabla u_n(x) \to \nabla u(x)\ \ {\rm a.e.} \ x\in B_R.
\end{eqnarray*}
Taking into account of the arbitrariness of $B_{ R}$, we conclude that
 \begin{eqnarray*}
 \nabla u_n(x)\to \nabla u(x) \ \ {\rm a.e.} \ x\in\mathbb{R}^3.
 \end{eqnarray*} We complete the proof. \hfill $\Box$

 Arguing as in \cite[Lemma 1.21]{1996Willem}, we   obtain the following lemma.
 \begin{lem}\label{lem3.15} Assume $\delta>0$ and $p\leqslant l<p^*$. Let $\{u_n\}\subset\W$ be a bounded sequence such that
\begin{eqnarray*}
\sup_{y\in\mathbb{R}^3}\int_{B(y,\delta)}|u_n|^l dx\rightarrow0.
\end{eqnarray*}
Then $u_n\rightarrow0$ in $L^{\beta}({\mathbb{R}^3})$ for any $p<\beta<p^*$.
\end{lem}

The next lemma is a key part of proving the existence of nontrivial solutions of \eqref{101}.
\begin{lem}\label{lem3.14}  If $\{u_n\}\subset \W$ is bounded  and satisfies
\begin{eqnarray}\label{3.39}
\mathcal{J}(u_n) \rightarrow c>0,  \ \ \  \ \ \mathcal{J}'(u_n)\rightarrow 0,
\end{eqnarray}
 then there is $\tilde{u}\in \W\setminus\{0\}$ such that $\mathcal{J}'(\tilde{u})=0$.  \end{lem}
 {\it Proof } \ Since $\{u_n\}$ is bounded in $\W$, there exist $\delta>0$ and $\mu\geqslant0$ such that
\begin{eqnarray*}
\sup_{y\in\mathbb{R}^3}\int_{B(y,\delta)}|u_n|^pdx\rightarrow \mu.
\end{eqnarray*}
 If $\mu=0$, then Lemma \ref{lem3.15} implies that
 \begin{eqnarray}
\int_{\mathbb{R}^3}|u_n|^\beta dx\rightarrow0, \ \ \ \forall \ p<\beta<p^*.\label{3.1}
\end{eqnarray}
 Thus
\begin{eqnarray}\label{3.2}
\int_{\mathbb{R}^3} \phi_{u_n}|u_n|^sdx\leqslant C \left(\int_{\mathbb{R}^3}|u_n|^{\frac{q^*s}{q^*-1}} dx\right)^{\frac{(q^*-1)q}{q^*(q-1)}}\rightarrow0.
\end{eqnarray}
By \eqref{3.39}   we get that
\begin{eqnarray*}\label{7.7}
\lambda\left(\frac{q-1}{qs}-\frac{1}{p}\right)\int_{\mathbb{R}^3}\phi_{u_n} |u_n|^sdx+\frac{r-p}{rp}\int_{\mathbb{R}^3} |u_n|^rdx=\mathcal{J}(u_n)-\frac{1}{p}\mathcal{J}'(u_n)[u_n]\rightarrow c.
\end{eqnarray*}
 From \eqref{3.1} and \eqref{3.2} we see that $c=0$ which contradicts  $c>0$. Thus $\mu>0$.   Then there exists $\{y_n\}\subset\mathbb{R}^3$ such that
\begin{eqnarray}\label{3.0}
\int_{B(y_n,\delta)}|u_n|^pdx \geqslant \frac \mu 2>0.
\end{eqnarray}
Let $\tilde{u}_n(x)=u_n(x+y_n)$. We deduce from  \eqref{3.0} that
\begin{eqnarray}\label{2.9}
\int_{B(0,\delta)}|\tilde{u}_n|^pdx \geqslant\frac \mu 2>0.
\end{eqnarray}
 It is obvious that $\{\tilde{u}_n\}$ is bounded in $\W$. Then, up to a subsequence, there exists $\tilde{u}\in\W$ such that
\begin{eqnarray}
  \left\{\begin{array}{ll}
\tilde{u}_n\rightharpoonup \tilde{u} &  \mbox{in} \ \W;\\
\tilde{u}_n\rightarrow \tilde{u} &  \mbox{in} \ L_{\rm loc}^l({\mathbb{R}^3})  \ \ \mbox{for} \ p\leqslant l<p^*;\\
\tilde{u}_n(x)\rightarrow \tilde{u}(x) &  \mbox{a.e.\ } \  {x\in\mathbb{R}^3}.
  \end{array}
\right.\label{2.2}
\end{eqnarray}
From \eqref{2.9} and \eqref{2.2}, we see that
\begin{eqnarray*}
\int_{B(0,\delta)}|\tilde{u}|^pdx\geqslant\frac \mu 2>0.
\end{eqnarray*}
Therefore $\tilde{u}\neq0$.
 Since $\tilde{u}_n=u_n(\cdot+y_n)\rightharpoonup \tilde{u}$ in $\W$ and $\mathcal{J}'(u_n)\rightarrow 0$, by Proposition \ref{Pro2}(ii), it can be verified that $\mathcal{J}'(\tilde{u}_n)\rightarrow 0$. Then by \eqref{2.2}, using Lemma \ref{lem1} and Proposition \ref{Pro2}(iv), from \cite[Proposition 5.4.7]{2013Willem}, we can obtain that $\mathcal{J}'(\tilde{u})=0$. The proof is complete. $\hfill\Box$

Now we verify that the functional $\mathcal{J}$ has the mountain pass geometry.
\begin{lem}\label{lem39} Let $\frac{pq(1+s)}{p(q-1)+q}<r<p^*$ and $\lambda>0$, or $p<r\leqslant\frac{pq(1+s)}{p(q-1)+q}$ and $\lambda>0$ be small. Then there exist $\rho>0$ and $\sigma>0$ such that
\begin{enumerate}
  \item[{\rm (i) }] there is $v\in \W$ so that $\|v\|> \rho$ and $\mathcal{J}(v)<0;$
     \item[{\rm (ii) }] $\mathcal{J}(u)\geqslant \sigma$ for any $\|u\|= \rho$.
  \end{enumerate}
   \end{lem}
{\it Proof} \ (i) Let $\frac{pq(1+s)}{p(q-1)+q}<r<p^*$ and $\lambda>0$.
 For $u\in \W\setminus\{0\}$ and $t>0$, take $u_t(x)=t^{\frac{p(q-1)+q}{p(1-q)+qs}}u(tx)$. Using Proposition \ref{Pro2}(ii), a direct computation gives that
  \begin{eqnarray}
  \left.
  \begin{array}{ll} \displaystyle
 \mathcal{J}(u_t) & = \displaystyle\frac{t^{\alpha_1-3}}{p}\int_{\mathbb{R}^3}|\nabla u|^pdx+\frac{t^{\alpha_2-3}}{p}\int_{\mathbb{R}^3}|u|^pdx+\frac{\lambda(q-1) t^{\alpha_1-3}}{qs}\int_{\mathbb{R}^3}\phi_u|u|^sdx \\[3mm]
&  \displaystyle \ \ \  -\frac{t^{\alpha_3-3}}{r}\int_{\mathbb{R}^3}|u|^rdx,
 \end{array}
\right.\label{7.4}
 \end{eqnarray}
where
\begin{eqnarray}\label{2.7}
\alpha_1=\frac{pq(s+1)}{p(1-q)+q s},\ \ \alpha_2=\frac{p[p(q-1)+q]}{p(1-q)+qs},\ \
 \alpha_3=\frac{r[p(q-1)+q]}{p(1-q)+q s}.
 \end{eqnarray}
We claim that
    \begin{eqnarray}\label{2.1}
    \alpha_3> \alpha_1> \alpha_2,\ \ \  \alpha_1-3>0.
       \end{eqnarray}
       To see this, we first notice that  $\max\left\{1,\frac{(q^*-1)p}{q^*}\right\}<s<\frac{(q^*-1)p^*}{q^*}$ implies that
         \begin{eqnarray}\label{5.1}
      p(1-q)+qs>0,\ \ \  s<\frac{p(4q-3)}{q(3-p)}.
      \end{eqnarray}
As $1<p<3$ and $1<q<3$, using $r>\frac{pq(1+s)}{p(q-1)+q}$ and \eqref{5.1}, we deduce that
  \begin{eqnarray}\label{5.3}
\alpha_3>\alpha_1>3.
  \end{eqnarray}
An explicit calculation gives
   \begin{eqnarray}\label{5.4}
  \alpha_1-\alpha_2=p>1.
   \end{eqnarray}
From \eqref{5.3} and \eqref{5.4}, then \eqref{2.1} holds. Then for $t_0>0$ sufficiently large, (i) follows by taking $v=u_{t_0}.$

 Let $p<r\leqslant\frac{pq(1+s)}{p(q-1)+q}$. For   $u\in \W\setminus\{0\}$ and $t>0$,  by Proposition \ref{Pro2}(ii), we have that
  \begin{eqnarray*}
 \mathcal{J}(t u) = \frac{t^p}{p}\int_{\mathbb{R}^3}(|\nabla u|^p+|u|^p)dx+\frac{\lambda(q-1) t^{\frac{qs}{q-1}}}{q s}\int_{\mathbb{R}^3}\phi_u|u|^sdx  -\frac{t^r}{r}\int_{\mathbb{R}^3}|u|^rdx.
 \end{eqnarray*}
 Choosing $t_1>0$ sufficiently large and $\lambda>0$ small enough,  take $v=t_1u$, then (i) easily follows, since $r>p$.

 (ii) Using Proposition \ref{Pro2}(i) and Sobolev inequality, it is easy to see that
\begin{eqnarray*}
 \mathcal{J}(u)\geqslant \frac{1}{p}\|u\|^p- C\|u\|^r,
\end{eqnarray*}
from which, in view of $r>p$, we infer that (ii) holds. The proof is complete.$\hfill\Box$

By Lemma \ref{lem39} and $\mathcal{J}(0)=0$, the mountain pass level can be  defined as
\begin{eqnarray}\label{7.5}
\bar c=\inf_{\gamma\in \Gamma} \max_{t\in [0, 1]} \mathcal{J}(\gamma(t))>0,
 \end{eqnarray}
 where
 \begin{eqnarray}\label{7.6}
\Gamma =\left\{\gamma\in C([0, 1], \W): \ \gamma(0)=0, \ \mathcal{J}(\gamma(1))<0\right\}.
  \end{eqnarray}
Observe that Lemma \ref{lem39} and the mountain pass theorem \cite{1973AR} can produce a Palais--Smale sequence for $ \mathcal{J}$. However, for $p<r<\frac{qs}{q-1}$, the boundedness of the Palais--Smale sequences for $ \mathcal{J}$ can not be proved directly. We first treat the case $\frac{pq(1+s)}{p(q-1)+q}<r<\frac{q s}{q-1}$.
 We require an auxiliary functional $\widetilde{\mathcal{J}}: \W\rightarrow\mathbb{R}$ defined as
\begin{eqnarray*}
\widetilde{\mathcal{J}}(u)&=&\frac{\alpha_1-3}{p}\int_{\mathbb{R}^3}|\nabla u|^pdx+\frac{\alpha_2-3}{p}\int_{\mathbb{R}^3}|u|^pdx+\frac{\lambda(q-1)(\alpha_1-3)}{qs}\int_{\mathbb{R}^3}\phi_u |u|^sdx\\
&&-\frac{\alpha_3-3}{r}\int_{\mathbb{R}^3}|u|^rdx,
\end{eqnarray*}
 where $\alpha_i\ (i=1,2,3)$ are given in \eqref{2.7}. Using the technique in \cite{1997LJ}, we construct a bounded Palais--Smale sequence for $ \mathcal{J}$.
\begin{lem}\label{lem3.12} Let $\frac{pq(1+s)}{p(q-1)+q}<r<p^*$. Then there exists a bounded sequence $\{w_n\}\subset\W$ such that
\begin{eqnarray}\label{325}
\mathcal{J}(w_n)\rightarrow \bar{c}, \ \ \ \mathcal{J}'(w_n)\rightarrow 0 \ \ \  \mbox{and} \ \ \ \widetilde{\mathcal{J}}(w_n) \rightarrow 0.
\end{eqnarray}
\end{lem}
{\it Proof} \ Consider the map $K:\mathbb{R}\times\W\rightarrow\W$ defined as
$$ K(\tau,u(x))=e^{\frac{p(q-1)+q}{p(1-q)+qs}\tau}u(e^\tau x).$$
As in \eqref{7.4}, an explicit computation gives that
\begin{eqnarray}\label{326}
    \left. \begin{array}{ll}
    (\mathcal{J}\circ K)(\tau,u)& \displaystyle =\frac{e^{(\alpha_1-3)\tau}}{p}\int_{\mathbb{R}^3}|\nabla u|^pdx+\frac{e^{(\alpha_2-3)\tau}}{p}\int_{\mathbb{R}^3}|u|^pdx\\[3mm]
      & \displaystyle \ \ \ +\frac{\lambda(q-1) e^{(\alpha_1-3)\tau}}{qs}\int_{\mathbb{R}^3}\phi_u|u|^sdx-\frac{e^{(\alpha_3-3)\tau}}{r}\int_{\mathbb{R}^3}|u|^rdx. \end{array}\right.
 \end{eqnarray}
Using Proposition \ref{lem2.4},
it can be verified that $\mathcal{J}\circ K$ is continuously Fr\'{e}chet-differentiable on $\mathbb{R}\times \W$.
For $u\in \W\setminus\{0\}$, from \eqref{326} and \eqref{2.1}, there exists $\tau_0>0$ sufficiently large such that $(\mathcal{J}\circ K)(\tau_0,u)<0$. Then, by $(\mathcal{J}\circ K)(0,0)=0$, the minimax level for $\mathcal{J}\circ K$ can be defined as
\begin{eqnarray}
\displaystyle \check{c}=\inf_{\check{\gamma}\in\check{\Gamma}}\sup_{t\in[0,1]}(\mathcal{J}\circ K)(\check{\gamma}(t)),\label{327}
\end{eqnarray}
 where $$\check{\Gamma}=\{\check{\gamma}\in C([0,1],\mathbb{R}\times \W):\check{\gamma}(0)=(0,0),\ (\mathcal{J}\circ K)(\check{\gamma}(1))<0\}.$$
 By the definitions of $\Gamma$ and $\bar{c}$ given in \eqref{7.6} and \eqref{7.5}, we know that $\{K\circ\check{\gamma}:\check{\gamma}\in\check{\Gamma}\}\subset\Gamma$ and $\{0\}\times\Gamma\subset\check{\Gamma}$, from which we infer that
 \begin{eqnarray}\label{328}
 \check{c}=\bar{c}>0.
 \end{eqnarray}
By the definition of $\bar{c}$ given in \eqref{7.5}, there exists a sequence $\{\gamma_n\}\subset\Gamma$ so that
 \begin{eqnarray}
\sup_{t\in[0,1]}(\mathcal{J}\circ K)(0,\gamma_n(t))=\max_{t\in[0,1]}\mathcal{J}(\gamma_n(t))\leqslant \bar{c}+\frac{1}{n}.\label{2.3}
\end{eqnarray}
 From \eqref{327}--\eqref{2.3}, \cite[Theorem 2.8]{1996Willem} produces the sequence $\{(\tau_n,u_n)\}\subset \mathbb{R}\times \W$ satisfying
\begin{eqnarray}
  \left\{\begin{array}{ll}
&(\mathcal{J}\circ K)(\tau_n,u_n)\rightarrow \bar{c};\\
&(\mathcal{J}\circ K)'(\tau_n,u_n)\rightarrow 0;\\
&\displaystyle\inf_{t\in[0,1]}(|\tau_n|^2 +\|u_n-\gamma_n(t)\|^2)^{\frac12}\rightarrow 0.
  \end{array}
\right.\label{330}
\end{eqnarray}
For any $(\zeta,\upsilon)\in \mathbb{R}\times \W$, we have
\begin{eqnarray}
(\mathcal{J}\circ K)'(\tau_n,u_n)[\zeta,\upsilon]=  \mathcal{J}'(K(\tau_n,u_n))[K(\tau_n,\upsilon)]+ \widetilde{\mathcal{J}} (K(\tau_n,u_n))\zeta.\label{332}
\end{eqnarray}
 Let $w_n=K(\tau_n,u_n)$, from \eqref{330} and \eqref{332}, we see that $\{w_n\}\subset\W$ satisfies
 \eqref{325}. Moreover, we infer that for $n$ sufficiently large,
\begin{eqnarray*}
\left.
  \begin{array}{ll}
 \bar{c}+1 &\geqslant \displaystyle \mathcal{J}(w_n)-\frac{1}{\alpha_3-3}\widetilde{\mathcal{J}}(w_n) \\[3mm]
&\displaystyle =\frac{\alpha_3-\alpha_1}{p(\alpha_3-3)}\int_{\mathbb{R}^3}|\nabla w_n|^pdx+\frac{\alpha_3-\alpha_2}{p(\alpha_3-3)}\int_{\mathbb{R}^3} |w_n|^pdx \\[3mm] &  \ \ \ \displaystyle+\frac{\lambda(q-1)(\alpha_3-\alpha_1)}{qs(\alpha_3-3)}\int_{\mathbb{R}^3}\phi_{w_n} |w_n|^sdx.\\
  \end{array}
\right.
\end{eqnarray*}
Combining this with \eqref{2.1}, we obtain the boundedness of $\{w_n\}$ in $\W$. The proof is complete.$\hfill\Box$\\
{\it Proof of Theorem \ref{thm3.3}.}\ By Lemma \ref{lem3.12} and Lemma \ref{lem3.14}, the conclusion immediately follows. The proof is complete.$\hfill\Box$

\

  Now we give the proof of Theorem \ref{thm3.1}. According to Theorem \ref{thm3.3}, we only need to consider the case $p<r\leqslant\frac{pq(1+s)}{p(q-1)+q}$. In this case, the approach applied in Lemma \ref{lem3.12} is
  not effective for obtaining the  boundedness of the Palais--Smale sequences for $\mathcal{J}$.
  Instead, we apply an idea in \cite{2006JS} to construct a
  cut-off functional $\mathcal{J}_M: \W\to \mathbb{R}$   defined by
 \begin{eqnarray*} \mathcal{J}_M(u)=\frac{1}{p}\|u\|^p +\frac{\lambda(q-1)}{qs} Z_M(u) \int_{\mathbb{R}^3}\phi_u|u|^sdx-\frac{1}{r}\int_{\mathbb{R}^3}|u|^{r}dx,\label{3.29}
 \end{eqnarray*}
  where
$Z_M(u)=h\left(\frac{\|u\|^2}{M^2}\right)$,
$h\in C_0^{\infty}(\mathbb{R},[0,1])$ satisfies $h(t)=1$ for $t\in[0,\frac{1}{2}]$ and $\hbox{supp}\ h\subset[0,1]$. The functional $\mathcal{J}_M$ is of $C^1$ and
 \begin{eqnarray*}
 \mathcal{J}'_M(u)[u]=\|u\|^p+\lambda\left(Z_M(u) +\frac{2(q-1)\|u\|^2}{qs M^2}h'\left(\frac{\|u\|^2}{M^2}\right)\right)\int_{\mathbb{R}^3}\phi_u|u|^s dx-\int_{\mathbb{R}^3}|u|^{r}dx.
 \end{eqnarray*}
 For any $M>0$, the functional $\mathcal{J}_M$ possesses the mountain pass geometry.
\begin{lem}\label{lem2.3}  Assume $p<r<p^*$. There exist $\bar{\rho}>0$ and $\bar{\sigma}>0$ such that
  \begin{enumerate}
  \item[{\rm (i) }] there is $\bar{v}\in \W$ such that $\|\bar{v}\|> \bar{\rho}$ and $\mathcal{J}_M(\bar{v})<0;$
     \item[{\rm (ii) }] $\mathcal{J}_M(u)\geqslant \bar{\sigma}$ for all $\|u\|= \bar{\rho}$.
  \end{enumerate}
   \end{lem}
{\it Proof} \  The proof can be proceeded  as in Lemma \ref{lem39} in a simpler way, so we omit it here.
 $\hfill\Box$

 From Lemma \ref{lem2.3} and $\mathcal{J}_M(0)=0$, we define the mountain pass level for $\mathcal{J}_M$ as
 $$\bar c_M=\inf_{\gamma\in \Gamma} \max_{t\in [0, 1]} \mathcal{J}_M(\gamma(t))>0,$$
 where
 $$\Gamma_M =\left\{\gamma\in C([0, 1], \W): \ \gamma(0)=0, \ \mathcal{J}_M(\gamma(1))<0\right\}.$$
 We prove that for  $M>0$ large enough and $\lambda>0$ small, the Palais--Smale sequences of  $\mathcal{J}_M$ are the Palais--Smale sequences of  $\mathcal{J}$.
 \begin{lem}\label{lem3.88} There exists $\lambda^*:=\lambda^*(M)>0$
 such that, for  $M>0$ large enough and $0<\lambda<\lambda^*$, if $\{u_n\}\subset\W$ satisfies
 \begin{eqnarray} \label{3033}
\mathcal{J}_M(u_n) \rightarrow \bar{c}_M, \ \ \ \ \ \   \mathcal{J}_M'(u_n)\rightarrow 0,
\end{eqnarray} then
 \begin{eqnarray}\label{3.21}
\limsup_{n\rightarrow\infty}\|u_n\|<\frac{M}{2}.
  \end{eqnarray}
 \end{lem}
{\it Proof} \ We first prove that for any $M>0$, the  sequence $\{u_n\}$ is bounded in $\W$. Assume that
\begin{eqnarray}\label{3.25}
 \|u_n\|\rightarrow\infty \ \ \ \mbox{as} \
n\rightarrow\infty.
 \end{eqnarray}
 Then, from the definition of $\mathcal{J}_M$, we see that for $n$ sufficiently large,
\begin{eqnarray*}
 \mathcal{J}_M(u_n)=\frac{1}{p}\int_{\mathbb{R}^3}\left(|\nabla u_n|^p+|u_n|^p\right)dx-\frac{1}{r}\int_{\mathbb{R}^3}|u_n|^{r}dx.
 \end{eqnarray*}
 By \eqref{3033} we obtain that as $n$ large enough,
 \begin{eqnarray*}
\bar{c}_M+1+\|u_n \| \geqslant \mathcal{J}_M(u_n)-\frac{1}{r}  \mathcal{J}_M'(u_n)[u_n]
= \left(\frac{1}{p}-\frac{1}{r}\right)\|u_n\|^p,
\end{eqnarray*}
which contradicts \eqref{3.25} as $1<p<r$.
  To end the proof, reasoning by contradiction, we assume that, up to a subsequence,
\begin{eqnarray}
\lim_{n\rightarrow\infty}\|u_n\|\geqslant \frac{M}{2}.\label{3.23}
 \end{eqnarray}
On the one hand, using \eqref{3.23} and \eqref{3033}, by the boundedness of $\{u_n\}$, we deduce that for $n$  large enough,
\begin{eqnarray}\label{3.88}
\left(\frac{1}{p}-\frac{1}{r}\right)\|u_n\|^p+\frac{1} {r}\mathcal{J}_M'(u_n)[u_n]\geqslant CM^p-M.
\end{eqnarray}
On the other hand, by \eqref{3033}, we conclude that
\begin{eqnarray}
\left.
  \begin{array}{ll} \displaystyle
& \displaystyle \left(\frac{1}{p}-\frac{1}{r}\right)\|u_n\|^p+\frac{1} {r}\mathcal{J}_M'(u_n)[u_n]\\ \leqslant & \displaystyle  \mathcal{J}_M(u_n)+\left|\frac{\lambda(q-1)}{qs}-\frac{\lambda}{r}\right|Z_M(u_n) \int_{\mathbb{R}^3}\phi_{u_n}|u_n|^sdx\\ &
  \displaystyle \ \ \ +\frac{\lambda[2(q-1)\|u_n\|^2]}{r q s M^2}\left|h'\left(\frac{\|u_n\|^2}{M^2}\right)\right|\int_{\mathbb{R}^3} \phi_{u_n}|u_n|^sdx.
  \end{array}
\right.\label{3.22}
\end{eqnarray}
We then estimate the right hand side of \eqref{3.22}. Note that if $\|u_n\|>M$, then $Z_M(u_n)=0$ and $h'\left(\frac{\|u_n\|^2}{M^2}\right)=0.$ Thus, by Proposition \ref{Pro2}(iii) and $h\in C_0^{\infty}(\mathbb{R},[0,1])$, we have that
\begin{eqnarray}\label{3.67}
Z_M(u_n)\int_{\mathbb{R}^3}\phi_{u_n}|u_n|^s dx\leqslant C\|u_n\|^{\frac{qs}{q-1}}Z_M(u_n)\leqslant CM^{\frac{qs}{q-1}},
\end{eqnarray}
\begin{eqnarray}\label{3.68}
\|u_n\|^2\left|h'\left(\frac{\|u_n\|^2}{M^2}\right)\right|\int_{\mathbb{R}^3}\phi_{u_n}|u_n|^sdx\leqslant C\|u_n\|^{\frac{qs}{q-1}+2}\left|h'\left(\frac{\|u_n\|^2}{M^2}\right)\right|\leqslant C_1M^{\frac{qs}{q-1}+2}.
\end{eqnarray}
To estimate  $\mathcal{J}_M(u_n)$, we first estimate $\bar{c}_M$.
  Actually, in view of the definitions of $\bar{c}_M$ and $\bar{v}$, arguing as in \eqref{3.67}, we deduce that
\begin{eqnarray*}\label{3.59}
\bar{c}_M&\leqslant&\max_{t\in[0,1]}\mathcal{J}_M(t\bar{v})\\
&\leqslant&\max_{t\in[0,1]}\left\{\frac{t^p}{p}\|\bar{v}\|^p-\frac{ t^r}{r}\int_{\mathbb{R}^3}|\bar{v}|^{r}dx\right\}+\max_{t\in[0,1]}\left\{ \frac{\lambda(q-1) t^{\frac{qs}{q-1}}Z_M(t\bar{v})}{qs}\int_{\mathbb{R}^3}\phi_{\bar{v}} |\bar{v}|^s dx\right\}\nonumber\\
&\leqslant& C_2+C_3\lambda M^{\frac{qs}{q-1}}.
\end{eqnarray*}
Combining  this with the first part of  \eqref{3033},
we infer that for $n$  sufficiently large,
 \begin{eqnarray}\label{3.81}
 \mathcal{J}_M(u_n)\leqslant3\bar{c}_M\leqslant3C+3C\lambda M^{\frac{qs}{q-1}}.
     \end{eqnarray}
 Inserting \eqref{3.67}--\eqref{3.81} into \eqref{3.22}, we have that for $n$ sufficiently large,
 \begin{eqnarray}\label{3.7}
\left(\frac{1}{p}-\frac{1}{r}\right)\|u_n\|^p+\frac{1} {r}\mathcal{J}_M'(u_n)[u_n]\leqslant C+ C\lambda M^{\frac{qs}{q-1}}.
\end{eqnarray}
By \eqref{3.88} and \eqref{3.7}, we conclude that
\begin{eqnarray}\label{3.27}
CM^p-M\leqslant  C+ C\lambda M^{\frac{qs}{q-1}},
\end{eqnarray}
where $C>0$ does not depend on $M$ and $\lambda$.
Since $p>1$, for any $M>0$ large enough and $0<\lambda<M^{-\frac{qs}{q-1}}$, \eqref{3.27} leads to a contradiction. The proof is complete. $\hfill\Box$\\
{\it Proof of Theorem \ref{thm3.1}.}\ From Lemma \ref{lem2.3} and the mountain pass theorem \cite{1973AR}, there exists a sequence $\{u_n\}\subset \W$ such that
 \begin{eqnarray} \label{8.0}
\mathcal{J}_M(u_n) \rightarrow \bar{c}_M>0, \ \ \ \   \ \ \   \mathcal{J}_M'(u_n)\rightarrow 0.
\end{eqnarray}
In view of Lemma \ref{lem3.88}, let $M>0$ be sufficiently large and $\lambda^*=\lambda^*(M)>0$ be small, then for any $0<\lambda<\lambda^*$, \eqref{3.21} follows.
 Together with the definition of $\mathcal{J}_M$, we deduce from \eqref{8.0} that
  \begin{eqnarray*}
\mathcal{J}(u_n) \rightarrow \bar{c}_M>0, \ \ \ \   \ \ \   \mathcal{J}'(u_n)\rightarrow 0.
\end{eqnarray*}
Then the conclusion immediately follows from Lemma \ref{lem3.14}. The proof is complete.
     $\hfill\Box$

\section*{Acknowledgment}
This work is supported by NSFC(12271373,12171326), KZ202010028048 and Natural Science Foundation of Sichuan, China(2023NSFSC1298). The authors would like to thank the referees  for carefully reading the manuscript and giving valuable comments to improve the exposition of the paper.

 \end{document}